\documentclass{amsart}

\usepackage{amsthm,amsfonts,amsmath,amssymb,mathrsfs,verbatim,cite}

\numberwithin{equation}{section}
\newtheorem{theorem}{Theorem}[section]
\newtheorem{corollary}[theorem]{Corollary}
\newtheorem{proposition}[theorem]{Proposition}

\newtheorem{lemma}[theorem]{Lemma}

\newtheorem{definition}[theorem]{Definition}
\newtheorem{example}[theorem]{Example}

\newcommand{\Nat}{\mathbb N}
\newcommand{\Z}{\mathbb Z}
\newcommand{\F}{\mathbb F}

\DeclareMathOperator{\Ass}{Ass}
\DeclareMathOperator{\Mass}{Mass}

\DeclareMathOperator{\End}{End}
\DeclareMathOperator{\Hom}{Hom}
\DeclareMathOperator{\ann}{ann}

\begin{document}

\title{Right Mori orders}

\author{N H Halimi}


\address{School of Mathematics and Physics, 
The University of Queensland, QLD 4072, Australia}
\email{n.halimi@uq.edu.au}

\begin{abstract}
In this paper we study \textit{right Mori Orders}, which are those prime Goldie rings that satisfy the 
ascending chain condition on regular integral right divisorial right ideals.
We will show that the class of right Mori orders is closed with respect to Morita-equivalence.
We also prove that each regular right divisorial right ideal of a right Mori order is contained in only finitely many right divisorial
completely prime right ideals. Moreover, we show that such right divisorial ideals can be represent as a finite intersection of
$\nu$-irreducible ideals of the form $aS:_rb$ for some regular $a,b\in S$.
\end{abstract}
\maketitle
\section{introduction}
 
It is well known that, if the dimension of a commutative Noetherian domain $A$ is greater than two, then 
its integral closure $\bar{A}$ is not necessarily Noetherian.
However, $\bar{A}$ is completely integrally closed such that the ascending chain 
condition holds on the set of integral divisorial ideals, i.e., $\bar{A}$ is a Krull domain.
A domain that satisfies the ascending chain condition on the set of integral divisorial ideals is called a Mori domain. 
This notion was introduced in the 1970s by J.~Querr\'e \cite{Q76} and has attracted a lot of attention since that time,
see for example the works of Barucci, Gabell, Houston and Lucas \cite{BS87, B83, B2000, L04, HL88}.
In \cite{L04}, Lucas extended the concept of Mori domain to rings with zero divisors.
Following the terminology of \cite{L04}, a ring is called a Mori ring if it satisfies the ascending chain
condition on the set of divisorial regular ideals. 

The purpose of this paper is to extend the notion of commutative Mori ring to prime Goldie rings.
We will refer to such extensions as right Mori orders.
We also investigate those properties of Mori domains that can be carried over to Mori orders.  
In particular we study the relation between completely prime right ideals and maximal elements of certain types of sets.
One of the outcomes is that each right divisorial ideal can be represent as a finite intersection $\nu$-irreducible ideals.
Moreover, we examine the right Mori property of pairs of local orders with the same prime ideals, and show that
if one of them is right Mori then so is the other.

Let $S$ be an order in a simple Artinian ring $Q$, i.e., $S$ is a prime Goldie ring with total quotient ring $Q$.
Let $U(Q)$ be the group of units in $Q$. A right $S$-submodule $I$ of $Q$ is called a \textit{right $S$-ideal}
if $I$ contains a regular element in $S$ and $uI\subseteq S$ for some $u\in U(Q)$. For any pair of subsets
$A$ and $B$ of $Q$, we use the notations
$(A:B)_r=\{q\in Q: Bq\subseteq A\}$ and $(A:B)_l=\{q\in Q:qB\subseteq A\}$.
If $I$ is a right $S$-ideal, then $(S:I)_l$ is a left $S$-ideal. 
We denote $I_{\nu}=(S:(S:I)_l)_r$. The set $I_{\nu}$ is called a \textit{right $\nu$-ideal} if $I_{\nu}=I$.
Similarly, for any left $S$-ideal $J$, we can define a left $S$-ideal $_{\nu}J$. An $S$-ideal $I$ is called a $\nu$-ideal if $I_{\nu}=I=_{\nu}I$.
An order $S$ in $Q$ is called a \textit{right Mori} order if the ascending chain condition holds for regular integral right $\nu$-ideals.
\medskip
  
In the first part of this paper we  define the notion of right Mori order and establish its basic properties.
In Theorem~\ref{finite}, it is proved that each regular right divisorial ideal of a right Mori order is contained in only finitely many
right divisorial completely prime ideals. This is an extension of Lucas' \cite[Theorem 2.18]{L04} to the non-commutative situation.

\medskip

Section~\ref{Sec3}, is concerned with pairs of orders with the same prime ideals. Theorem~\ref{pair} shows that if one of them is right Mori
then so is the other. Moreover, it is proved that a non-commutative Krull order in the sense of Marubayashi is a Mori order. 

\medskip
In section~\ref{Sec4} we focus on $\nu$-irreducible ideals and prove that in a right Mori order $S$, each right divisorial ideal
can be represented as a finite intersection of $\nu$-irreducible ideal of the form $aS:_rb$ for some regular elements $a, b\in S$.

\medskip

In the final section, using an example of Cohn and Schofield from \cite{CS85}, we show that a right Mori order is not necessarily left Mori.

\section{definitions and basic properties}\label{Sec3}
Unless stated otherwise, in this paper $S$ will denote an order in a simple Artinian ring $Q$, i.e., 
$S$ is a prime Goldie ring with total quotient ring $Q$.
A proper right ideal $P$ of a ring $R$ is called a completely prime right ideal if 
$aP\subseteq P$ and $ab\in P$ implies that $a\in P$ or $b\in P$ for all $a,b\in R$.
Two sided completely prime ideals are examples of completely prime right ideals. 
For more details on completely prime right ideals, see \cite{MR10}.   
\begin{proposition}\label{M}
Let $S$ be an order and $M$ a right ideal of $S$ such that $M$ is maximal with respect to being right divisorial. 
Then $M$ is a completely prime right ideal. 
\end{proposition}
\begin{proof}
Proceeding by contradiction we assume there exist $a, b\in S-M$ such that $ab\in M$ and $aM\subseteq M$. By maximality of $M$
we have $(M+bS)^{\nu}=S$. By our assumption we have $a(M+bS)\subseteq M$, so that $a(M+bS)^{\nu}=(a(M+bS))^{\nu}\subseteq M$.
Hence $a\in aS=a(M+bS)^{\nu}=(a(M+bS))^{\nu}\subseteq M$, which shows that $a\in M$, a contradiction.
\end{proof}

\begin{proposition}
Let $S$ be an order which contains at least one proper right $\tau$-ideal. Then:
\begin{itemize}
\item[(1)]
Maximal right $\tau$-ideals do exist and are completely prime right ideals.
\item[(2)]
Each proper right $\tau$-ideal is contained in a maximal right $\tau$-ideal.

\end{itemize}
\end{proposition}
\begin{proof}
(1). By assumption the set of all proper right $\tau$-ideals is non-empty. 
Now let $\{J_{\alpha}\}$ be a chain of right $\tau$-ideals and $I$ be a finitely generated right ideal of $S$
such that $I\subseteq \cup J_{\alpha}$. Then $I\subseteq J_{\alpha}$ for $\alpha$ large enough. Since
$J_{\alpha}$ is a right $\tau$-ideal we have $I^{\nu}=I^{\tau}\subseteq J_{\alpha}$. Thus $I^{\nu}\subseteq \cup J_{\alpha}$,
and hence $\cup J_{\alpha}$ is a right $\tau$-ideal. 
By Zorn's lemma the set of all proper right $\tau$-ideals has a maximal element $P$. Similar to the proof of Proposition~\ref{M},
it follows that $P$ is a completely prime right ideal.

(2). The second claim follows from (1).
\end{proof}
Recall that a right ideal $I$ of a ring $R$ is called regular if and only if $I$ contains a regular element of $R$.
\begin{theorem}\label{chain}
For any prime Goldie ring $S$ the following conditions are equivalent:
\begin{itemize}
\item[(1)]
The ascending chain condition (ACC) holds for regular integral right divisorial ideals.
\item[(2)]
For any regular integral (resp. fractional) right ideal $I$ of $S$,
there exists a finitely generated integral (resp. fractional) right ideal $J\subseteq I$ such that  $I^{\nu}=J^{\nu}$.
\end{itemize}
\end{theorem}

\begin{proof}
$(1)\Rightarrow (2)$.
It is enough to assume that $I$ is an integral right ideal of $S$, because if $I$ is fractional right ideal, then $uI\subseteq S$
for some $u\in U(Q)$. Thus let $I$ be an integral right ideal of $S$ and $\Gamma$ be the set of all $J_{\alpha}^{\nu}$, 
where $J_{\alpha}$ is finitely generated right ideal of $S$ with $J_{\alpha}\subseteq I$. 
The set $\Gamma$ has a maximal element $J^{\nu}$, 
because ACC holds for integral right divisorial ideals of $S$. If $J^{\nu}\subset I^{\nu}$,
then there exists an element $b\in I-J^{\nu}$. Put $J'=J+bS$. Then $J'\subseteq I$ and $J'$ is a finitely generated right ideal.
Thus $J'^{\nu}\in \Gamma$ with $J^{\nu}\subset J'^{\nu}$, in contradiction with the maximality of $J^{\nu}$.
Hence $J^{\nu}=I^{\nu}$.

$(2)\Rightarrow (1)$. Let $\{I_n\}$ be an ascending chain of regular integral right divisorial ideals of $S$.
Put $I=\cup I_n$. By assumption there exists a finitely generated right ideal $J\subseteq I$ with $J^{\nu}=I^{\nu}$. 
Since $I$ is the union of the chain $\{I_n\}$ and $J$ is finitely generated, we have $J\subseteq I_m$ for some positive $m$.
Since $I_m$ is right divisorial, we have $I_m\subseteq I\subseteq I^{\nu}=J^{\nu}\subseteq I_m^{\nu}=I_m$.
Hence $I_m=I$ and the chain $\{I_n\}$ stabilizes. 
\end{proof}
\begin{definition}
A prime Goldie ring satisfying the conditions of Theorem \ref{chain} is called a \textit{right Mori} order. 
\end{definition}
A left Mori order is defined similarly. A Mori order is an order which is both right and a left Mori order.
Commutative Mori domains,  commutative Krull domains and right Noetherian rings are examples of right Mori orders.
We will later prove that a non-commutative Krull order in the sense of Marubayashi is also a Mori order.

\begin{lemma}\label{DCC}
Let $S$ be a right Mori order. Then any descending chain of right divisorial regular ideals with regular intersection stabilizes. 
\end{lemma}
\begin{proof} 
Let $\{I_n\}$ be a descending chain of regular right divisorial ideals of $S$ with regular intersection. 
Put $I=\cap I_n$. Then $\{(S:I_n)_r\}$ is an ascending chain of regular right divisorial fractional ideals of $S$ such that $(S:I_n)_r\subseteq (S:I)_r$
for all $n$.
By assumption there exists a finitely generated regular fractional right ideal $J_n\subseteq (S:I_n)_r$ of $S$ such that $J_n^{\nu}=(S:I_n)_r$.
Hence we have an ascending chain of finitely  generated regular fractional right ideal $\{J_n\}$.
For $J=\cup J_n$ there exists a finitely generated fractional right ideal $K\subseteq J$ such that $K^{\nu}=J^{\nu}$. 
By construction of $J$, we have $K\subseteq J_m$ for some $m$.
Thus $K\subseteq J_m\subseteq J\subseteq J^{\nu}=K^{\nu}\subseteq J_m^{\nu}=J_m$ so that $J_m=J$. Therefore, 
$(S:I_m)_r=(S:I_i)_r$ for all $i\geq m$. From the fact that $I_n$ right divisorial, we can conclude that the 
descending chain of regular right divisorial ideals $\{I_n\}$ must stabilizes at $I_m$.
\end{proof}
In the commutative setting the converse of the above lemma is also true but
in the non-commutative case we were unable to prove or disprove such a converse result.
 
An order $S$ is called local if the Jacobson radical $J(S)$ is the only maximal right (left) ideal of $S$.
A local Bezout order in a simple Artinian ring $Q$ is called a Dubrovin valuation ring.
A Dubrovin valuation ring $S$ is called discrete if it is not Artinian and $\cap_{n=1}^{\infty} J(S)^n=0$,
see \cite[Definition 1.16]{AD90}. 
\begin{lemma}\label{D}
Let $S$ be a Mori order in a simple Artinian ring $Q$. If $S$ is a Dubrovin valuation ring
then $S$ is a rank one discrete valuation ring.
\end{lemma}
\begin{proof}
It is easy to show that $S$ is a principal ideal ring. Thus every ideal is  divisorial, so that
$S$ is Noetherian. By \cite[Proposition 5.16]{MMU97} $S$ is a local Dedekind ring and by \cite[Theorem 2.7]{K72}
$S$ is discrete. 
\end{proof}
A ring $S$ is called right \textit{quasi-coherent}
if the intersection of finitely many principal right ideals is a finitely generated right
ideal and for any $a\in S$ the right ideal $\ann_r(a)=\{x\in S: ax=0\}$ is finitely generated. 
By \cite [4.60 Corollary] {L99}, any right coherent ring is right quasi-coherent. Thus right Noetherian rings,
Dubrovin valuation rings, right Bezout, right Pr\"ufer rings and right semiherditary rings are all examples of right quasi-coherent rings.

\begin{proposition}
Let $S$ be an order in a simple Artinian ring $Q$ such that $\ann_r(a)$ is a finitely generated right ideal for any $a\in S$.
Then:
\begin{itemize}
\item[(1)]
If each right divisorial $I$ is finitely generated right ideal, then $S$ is right pseudo-coherent;
\item[(2)]
If $S$ is right Mori, then each right divisorial $I$ is finitely generated right ideal if and only if
$S$ is right pseudo-coherent.
\end{itemize}
\end{proposition}
\begin{proof}
(1).
$S$ is right pseudo-coherent, because the intersection of finitely many right principal ideals is a right divisorial ideal.

(2).
Let $S$ be a right pseudo-coherent ring and $I$ be a right divisorial right ideal of $S$. Since $S$ is a right Mori by 
Proposition~\ref{P2.8}, $I=\cap_{i=1}^n u_iS$ for some $u_i\in U(Q)$. Hence $I$ is a finitely generated right ideal of $S$. 
\end{proof}

The following theorem is an extension of \cite[Theorem 2.18]{L04} by Lucas to the non-commutative situation. 
\begin{theorem}\label{finite}
Let $S$ be a right Mori order. Then each regular right divisorial ideal is contained in only finitely many
right divisorial completely prime ideals.
In particular, each regular right ideal is contained in at most finitely many maximal right $\tau$-ideals. 
\end{theorem}
\begin{proof}
Let $\{P_i\}$ be a family of right divisorial completely prime ideals which contain of regular right divisorial ideal $I$. 
We proceed by contradiction and assume that $\{P_i\}$ is infinite. Without loss generality we can assume that
$\{P_i\}$ is countable. Since $S$ is right Mori, every chain in $\{P_i\}$ is finite.
Thus we can assume that for each $i\neq j, ~P_i$ and $P_j$ are not comparable. Put $I_n=\cap_{i=1}^nP_i$.
By assumption each $P_i$ is completely prime. Thus $I_{n+1}$ is a proper sub-ideal of $I_n$ for all $n\geq 1$.
Since the intersection of right divisorial ideals is again right divisorial, each $I_n$ is a right divisorial.
Hence $\{I_n\}$ is an infinite descending chain of regular right divisorial with regular intersection,
a contradiction with Lemma~\ref{DCC}. Since in a right Mori order a right $\tau$-ideal is a right divisorial
ideal, the rest of the proof follows from Proposition~\ref{M}.
\end{proof}

An $S$-ideal $I$ is called divisorial if $I=I^{\nu}={^{\nu}I}$. 
An integral ideal $I$ is called maximal divisorial if $I$ is maximal with respect to being
divisorial. By Proposition~\ref{M}, each maximal divisorial ideal is a completely prime ideal.

The following proposition can be applied to the decomposition of a Mori order, 
see for example, \cite[Proposition 2.2]{BS87} for the commutative case. 
\begin{proposition}
Let $S$ be a Mori order in a simple Artinian ring $Q$ such that $S$ is localizable at every $P\in D_m(S)$, where
$D_m(S)$ denotes the set of all maximal divisorial ideals.
Let $I$ be a right $S$-ideal. Then $IS_P=S_P$ for all but finitely many $P\in D_m(S)$.
\end{proposition}
\begin{proof}
 Let $I$ be a regular integral right ideal of $S$. If $I^{\nu}=S$, then $I\nsubseteq P$ for all $P\in D_m(S)$.
Thus for every $P\in D_m(S)$ there exists an element $a\in I-P$ such that $a+P$ is a regular element of $S/P$.
Therefore, $IS_P=S_P$. 
If $I^{\nu}\neq S$, then by Theorem \ref{finite}, $I$ contains only finitely many primes, say $P_1,\dots, P_n\in D_m(S)$.
Therefore, for each $P\in D_m(S)-\{P_1,\dots,P_n\}$ there exists $a\in I-P$ such that $a+P$ is regular
in $S/P$. Hence $IS_P=S_P$ for all $P\in D_m(S)-\{P_1,\dots,P_n\}$. Now let $I$ be a right $S$-ideal. Then
$qI\subseteq S$ for some regular element $q\in Q$. By the above, $qIS_P=S_P$ for all but finitely many $P\in D_m(S)$.
Since $Q$ is a simple Artinian ring and $q$ is regular, we have $q\in U(Q)$. Thus $IS_P=q^{-1}S_P$ 
for all but finitely many $P\in D_m(S)$. Since $S$ is an order in $Q$ we have $q^{-1}=b^{-1}a$ for
some $a, b \in S$ with $b$ regular. Again by the above, we have $aS_P=S_P$ for all but finitely many $P\in D_m(S)$.
Hence $b^{-1}aS_P=b^{-1}S_P$ for all but finitely many $P\in D_m(S)$. But we know that $b^{-1}S_P=S_P$
for all but finitely many $P\in D_m(S)$.
Therefore, $q^{-1}S_P=b^{-1}aS_P=S_P$ for all but finitely many $P\in D_m(S)$. 
\end{proof}

\section{Local orders with the same prime ideals}\label{Sec4}
Given a pair of commutative local domains with the same prime ideals, it is well known that if one of them is a Mori domain 
then so is the other. In this section, we first extend the above property to the non-commutative situation. Then
we focus on studying the right Mori property of a family of overrings of a ring with finite character. By applying 
Theorem~\ref{ACC} below  we conclude that the non-commutative Krull order in the sense of Marubayashi is a Mori order.

We will start this section with a basic lemma. 
\begin{lemma}\label{div}
Let $S$ be a local order with maximal ideal $M\neq (0)$. Then:
\begin{itemize}
\item[(1)]
$M$ is right divisorial if and only if $S\subset (S:M)_l$.
\item[(2)]
$M$ is a principal ideal if and only if $O_l(M)\subset (S:M)_l$.
\item[(3)]
$M$ is a non-principal right divisorial ideal if and only if $S\subset O_l(M)$. 
In this case $(S:M)_l=O_l(M)$.

\end{itemize} 
\end{lemma}
\begin{proof}

(1).
Since $M$ is a two sided ideal we always have $S\subseteq (S:M)_l$. Now if $S=(S:M)_l$,
then $M^{\nu}=(S:(S:M)_l)_r=(S:S)_r=S\neq M$, a contradiction with right divisoriality of $M$. 

Conversely, let $S\subset (S:M)_l$ and $x\in (S:(S:M)_l)_r$. Then $(S:M)_lx\subseteq S$ and so $x\in S$.
If $x\notin M$, then $x$ is a unit in $S$. For any $r\in (S:M)_l$, we have $rx\in S$.
Thus $r\in Sx^{-1}=S$ and so $(S:M)_l\subseteq S$, which is a contradiction.

(2).
Let $M=aS=Sa$ for some regular element $a\in M$. Then $M$  is a right divisorial ideal, and by part (1), $S\subset (S:M)_l$.
Now let $x\in O_l(M)$. Then $xaS\subseteq Sa$. Therefore, $xa=ra$ for some $r\in S$, which shows that $x=r\in S$
and $O_l(M)\subseteq S\subset (S:M)_l$. Conversely, let $x\in (S:M)_l-O_l(M)$. Then $xM\subseteq S$. Now if
$xM\neq S$, then by locality of $S$, we have $xM\subseteq M$, which is a contradiction. Thus $xM=S$ and so
$M=x^{-1}S$.

(3).
By part (2) $M$ is not principal if and only if $O_l(M)=(S:M)_l$. By part (1) $M$ is divisorial if and only
if $S\subset (S:M)_l$. Therefore, $M$ is a non-principal right divisorial ideal if and only if $S\subset O_l(M)$.
\end{proof}
For an order $T$ in a simple Artinian ring $Q$ we always assume that an overing $T$ is contained in $Q$.
\begin{lemma}\label{s}
Let $T\subset S$ be orders. If $J(T)$ is an ideal of $S$, then $J(T)\subseteq J(S)$. 
In particular if $T$ is a local ring with  maximal ideal $M$, then $M\subseteq J(S)$. 
\end{lemma}
\begin{proof}
The proof is the same as in the commutative case, see the proof of \cite[Lemma 3.6]{AD80}.
\end{proof}
The next two propositions are non-commutative versions of \cite[Proposition 2.4 and 2.6]{B83}  
and describe the sets of right divisorial ideals of $T$ and $S$.
\begin{proposition}\label{NP}
Let $T\subset S$ be local orders in a simple Artinian ring $Q$ with the same prime ideals.
Then each non-principal right divisorial ideal of $S$ is a right divisorial ideal of $T$.
\end{proposition}
\begin{proof}
Let $I$ be a non-principal right divisorial ideal of $S$. Then $I=\cap \{xS: I\subseteq xS,~ x\in U(Q)\}$.
Since $I$ is not principal, the inclusion $I\subseteq xS$ implies that $I\subset xS$.
Because $M$ is a maximal ideal of $S$,
the inclusion $I\subset xS$ implies that $I\subseteq xM$. Thus $I=\cap \{xM: I\subseteq xM,~ x\in U(Q)\}$. 
By Lemma \ref{div} part (3), $M$ is a right divisorial ideal of $T$, because $T\subset S\subseteq O_l(M)$.
Hence $xM$ is a right divisorial ideal of $T$ for all $x\in U(Q)$.
Since the intersection of right divisorial ideals is again right divisorial, $I$ is a right divisorial ideal of $T$.
\end{proof}
\begin{proposition}\label{xM}
Let $T\subset S$ be local orders in a simple Artinian ring $Q$ with the same prime ideals.
Then each non-principal right divisorial ideal $I$ of $T$ is a right fractional ideal of $S$ of at least one of the following types:
\begin{itemize}
\item[(1)]
$I=xM$ with $0\neq x\in Q$;

\item[(2)]
$I$ is a right divisorial ideal of $S$.
\end{itemize}
\end{proposition}

\begin{proof}
Let $I$ be a non-principal right divisorial ideal of $T$. We first show that $I$ is a right ideal of $S$,
that is $IS\subseteq I$. For this it is enough to show that $IS\subseteq (T:(T:I)_l)_r$. Now let $y\in (T:I)_l$.
Then $yI\subseteq T$. Since $I$ is not right principal, $yI$ is contained a maximal 
right ideal of $N$ of $T$. The ring $T$ is local so that $N=M$ and $yI\subseteq M$.
Thus $yIS\subseteq MS=M\subset T$. Hence $(T:I)_lIS\subseteq T$ so that $IS\subseteq (T:(T:I)_l)_r=I$.
If $I\neq xM$ for all $x\in Q$, we will show that $I$ is a right divisorial ideal of $S$. To see this 
it is enough to show that $(T:I)_la\subseteq T$ for all $a\in (S:(S:I)_l)_r$. Let $a\in (S:(S:I)_l)_r$.
Then from $(T:I)_l\subseteq (S:I)_l$, we conclude that $(T:I)_la\subseteq (S:I)_la\subseteq S$. Since 
$M$ is a right divisorial ideal of $T$ and $(T:M)_l=O_l(M)$, we have $M=(T:S)_r=(T:O_l(M))_r$.

Furthermore, since $I=(T:(T:I)_l)_r\neq aM$ and $aM=a(T:S)_r=(T:Sa^{-1})$, we have $(T:(T:I)_l)_r\neq (T:Sa^{-1})_r$.
Thus $(T:I)_l\neq Sa^{-1}$ so that $(T:I)_la\subset S$. By locality of $S$ we have $(T:I)_la\subseteq M\subset T$, as desired. 
\end{proof}
The following theorem is an extension to the non-commutative case of \cite[Theorem 3.2]{B83}.
\begin{theorem}\label{pair}
Let $T\subset S$ be local orders in a simple Artinian ring $Q$ with the same prime ideals.
Then $T$ is a right Mori order if and only if $S$ is a right Mori order.
\end{theorem}
\begin{proof}
Let $T$ be a right Mori order. By Proposition \ref{NP}, to prove that $S$ is a right Mori order 
it is enough to show that ACC holds for regular right principal integral ideals of $S$.
Let $s_1S\subset s_2S\subset \dots$ with $s_n\in M$ an increasing sequence of regular integral principal ideals of $S$.
Since $M$ is a right divisorial ideal of $T$ and $O_l(M)=(T:M)_l$, the left order $O_l(M)$ is a right divisorial $T$-ideal.
For any $n\in \Nat,~ s_nT\subseteq MT=M\subset T$ is an integral right divisorial ideal of $T$.
Therefore, there exists $n_0\in \Nat$ such that $s_nT=s_{n+1}T$ for all $n\geq n_0$.
Since $s_n$ is a regular element of $S$ for all $n\in \Nat$, the inverse of $s_n$ exists and $s_{n+1}^{-1}s_n\in U(T)$ for all
$n\geq n_0$. Therefore, $s_{n+1}^{-1}s_n\in U(S)$ and $s_nS=s_{n+1}S$ for all $n\geq n_0$.
Conversely, let $S$ be a Mori order. By Proposition~\ref{xM}, it is enough that to prove that ACC holds for 
regular principal right ideals and regular right ideal of the $xM$, where $0\neq x\in T$. 
Let $x_1T\subset x_2T\subset\dots $. Then $x_1S\subset x_2S\subset\dots $.
By assumption there exists an element $n_0\in \Nat$ such that $x_{n_0}S=x_nS$ for all $n\geq 0$.
Thus $x_{n+1}^{-1}x_n\in U(S)\cap T=U(T)$ for all $n\geq n_0$,
which shows that $x_{n+1}^{-1}x_nT=T$ and $x_{n+1}T=x_nT$ for all $n\geq n_0$.

Now assume that $x_1M\subset x_2M\subset\dots$ is an increasing chain of right ideals of $T$.
For each $n,~ x_nM$ is a right ideal of $S$, and $x_nM$ is a right divisorial ideal of $S$ if and only if
$M$ is a right divisorial ideal of $S$. By parts (2) and (3) of Lemma~\ref{chain}, $x_nM$ is a right divisorial ideal of $S$
for all $n\in \Nat$. Now, since $S$ is a right Mori order, there exists $n_0\in \Nat$ such that $x_nM=x_{n+1}M$ for all $n\geq n_0$. 
\end{proof}

The following is the definition of a non-commutative Krull ring in the sense of Marubayashi \cite{M75}. 
In the remainder of this section we will prove that non-commutative Krull rings are examples of Mori orders.

\begin{definition}
An order $S$ in a simple Artinian ring of $Q$ is called a Krull order if there exist 
families $\{R_i\}_{i\in \alpha}$ and $\{S_j\}_{j\in \beta}$ of essential overrings of $S$ such that:
\begin{itemize}
\item[(K1)]
$S=(\cap_{i\in \alpha} R_i)\cap (\cap_{j\in \beta} S_j)$;
\item[(K2)]
Each $R_i$ is a non-commutative discrete valuation ring, $S_j$ is a simple Noetherian ring and $|\beta|< \infty$;
\item[(K3)]
For every regular element $c\in S$ we have $cR_i\neq R_i$ ($R_ic\neq R_i$) for only finitely many
$i\in \alpha$.
\end{itemize}
$S$ is called bounded if $\beta=\emptyset$.
\end{definition}
A family of overings $\{S_i\}_{i\in \alpha}$ of $S$ with $S=\cap_{i\in \alpha}S_i$ is called of \textit{finite character}
if every non-zero non-unit element of $S$ is a non-unit of finitely many $S_i$.
\begin{lemma}\label{cap}
Let $\{S_i\}$ with $S=\cap_{i\in \alpha} S_i$ be a family of overrings of finite character, and let $I$ and $J$ be
right and left fractional ideals of $S$ respectively. Then $(S:I)_l=\cap_{i\in \alpha} (S_i :IS_i)_l$ and 
$(S:J)_r=\cap_{i\in \alpha} (S_i :S_iJ)_r$, so that $I^{\nu}=\cap _{i\in \alpha}(S_i :S_i(S:I)_l)_r$.  
\end{lemma}
\begin{proof}
Let $x\in (S:I)_l$. Then $xI\subseteq S$ and so $xIS_i\subseteq S_i$ for all $i\in \alpha$.
Hence $(S:I)_l\subseteq \cap_{i\in I} (S_i :IS)_l$. Conversely, let $x\in \cap_{i\in \alpha} (S_i :IS_i)_l$.
Then $xIS_i\subseteq S_i$ for all $i\in \alpha$. Now since $I\subseteq IS_i$ for all $i\in \alpha$, we have 
$xI\subseteq xIS_i\subseteq S_i$ for all $i\in \alpha$. Hence $xI\subseteq \cap_{i\in \alpha} S_i=S$,
which shows that $\cap_{i\in \alpha} (S_i :IS_i)_l\subseteq (S:I)_l$. Thus $(S:I)_l=\cap_{i\in I} (S_i :IS_i)_l$. 

In a similar manner one can prove that $(S:J)_r=\cap_{i\in \alpha} (S_i :S_iJ)_r$. In particular, $(S:I)_l$
is a left fractional ideal of $S$. Hence $I^{\nu}=(S:(S:I)_l)r=\cap _{i\in \alpha}(S_i :S_i(S:I)_l)_r$.
\end{proof}
\begin{lemma}\label{ex}
Let $\{S_i\}$ with $S=\cap_{i\in \alpha} S_i$ be a family of finite character and let $I\subset J$ be two right $\nu$-ideals of $S$.
Then there exists $i\in \alpha$ such that $(S_i:S_i(S:I)_l)_r\subset (S_i:S_i(S:J)_l)_r$.
\end{lemma}
\begin{proof}
Since $I\subset J$ are $\nu$-ideals, we have $(S:J)_l\subset (S:I)_l$.
Thus $S_i(S:J)_l\subset S_i(S:I)_l$ for all $i\in \alpha$, so that $(S_i: S_i(S:I)_l)_r\subseteq (S_i: S_i(S:J)_l)_r$.

If for all $i\in\alpha$ we have $(S_i: S_i(S:I)_l)_r=(S_i: S_i(S:J)_l)_r$, then by Lemma~\ref{cap}, 
$I=I^{\nu}=\cap _{i\in \alpha}(S_i: S_i(S:I)_l)_r= \cap _{i\in \alpha}(S_i: S_i(S:J)_l)_r=J=J^{\nu}$, a contradiction.
\end{proof}
 
\begin{theorem}\label{ACC}
Assume that $\{S_i\}$ with $S=\cap_{i\in \alpha} S_i$ is a family of finite character.
If $S$ does not satisfy ACC on the set of regular right divisorial ideals, then neither does at least one of the $S_i$.
\end{theorem}
\begin{proof}
Since $S$ does not satisfy ACC on the set of regular right divisorial ideals, there exists
an infinite ascending chain  $I_1\subset I_2\subset \dots $ of regular right divisorial ideals of $S$. 
From $I_j\subseteq S$, we can conclude that
$(S_i:S_i(S:I_j)_l)_r\subseteq S_i$. Hence $(S_i:S_i(S:I_j)_l)_r$ is a right regular divisorial ideal of $S_i$. 
The family $\{S_i\}$ is of finite character.
Hence there exists only finite number of $S_i$ for which $(S_i :S_i(S:I_1)_l)_r\subset S_i$.
Without loss of generality we can assume that $F=\{S_1\dots S_n\}$ is the set of all $S_i$ such that $(S_i :S_i(S:I_1)_l)_r\subset S_i$. 
For any $k>1$, we have $(S_i :S_i(S:I_1)_l)_r\subseteq (S_i :S_i(S:I_k)_l)_r$. 
Thus $(S_i :S_i(S:I_k)_l)_r\subset S_i$ implies that $S_i\in F$. Now, by Lemma~\ref{ex}, for each $n$
there exists $S_i\in F$ such that $(S_i:S_i(S:I_n)_l)_r\subset (S_i:S_i(S:I_{n+1})_l)_r$. Since $I_1\subset I_2\subset \dots $
is an infinite ascending chain and $F$ is a finite set, there exists $S_i\in F$
such that $(S_i:S_i(S:I_1)_l)_r\subset (S_i:S_i(S:I_2)_l)_r\subset \dots $. 
\end{proof}
\begin{corollary}
Assume that $\{S_i\}$ with $S=\cap_{i\in \alpha} S_i$ is a family of finite character.
If each $S_i$ is a right Mori order then $S$ is a right Mori order. In particular, the intersection
of finitely-many right Mori orders is again a right Mori order.
\end{corollary}

\begin{corollary}
Non-commutative Krull orders in the sense of Marubayashi are Mori orders.
\end{corollary}

\section{$\nu$-irreducible ideals}\label{Sec5}

For any subset $A$ of $S$ and $x\in S$ we define $A:_rx=\{s\in S: xs\in A\}$. The set $A:_lx$
is defined similarly. 
A completely prime right ideal $P$ of a ring $S$ is \textit{associated} to a right ideal $I$ of $S$
if $P=I:_rx$ for some element $x\in S-I$. 

\begin{lemma}\label{MM}
Let $S$ be a right Mori order and $I$ a right divisorial right ideal 
of $S$. Let $\Gamma= \{I:_ra : a\in S-I\}$. 
Then $\Gamma$ has a maximal element and is a completely prime right ideal.
\end{lemma}
\begin{proof}
Since the set $\Gamma$ has the ACC property, $\Gamma$ has a maximal element $M$.
Let $xy\in M$ with $xM\subseteq M$. 
Proceeding by contradiction assume that $x, y\notin M$.
Let $a\in S-I$ such that $M=I:_ra$. Since $x\notin M$, we have $ax\in S-I$.
Thus $I:_rax\in \Gamma$. Since $xM\subseteq M=I:_ra$, we have $axm\in I$ for any $m\in M$. 
Hence $M\subseteq I:_rax$. Since $y\notin M$ and $y\in I:_rax$. Together with $M\subseteq I:_rax$, 
this contradicts the maximality of $M$.  
\end{proof}

Following \cite{HL88} we denote by $\Ass_r(I)$ the set of all completely prime right ideals of the form
$I:_ra, ~a\in S-I$ and $\Mass_r(I)$ the set of all maximal element of $\Ass_r(I)$.
\begin{corollary}
Let $S$ be a right Mori order and $I$ be a right divisorial right ideal of $S$.
Then the set of maximal elements of $\Gamma= \{I:_ra : a\in S-I\}$ is equal to $\Mass_r(I)$.
\end{corollary}
\begin{proof}
By Lemma~\ref{M}, the set of maximal elements of $\Gamma$ is a subset of $\Mass_r(I)$.
Conversely, let $P\in \Mass_r(I)$. Then $P\in \Gamma$. Assume that $P$ is not a maximal element of $\Gamma$.
Then there exists an $a_1\notin I$ such that $P\subset I:_ra_1$. The right ideal $I:_ra_1$ is not 
a maximal element of $\Gamma$, because otherwise, by Lemma~\ref{M}, the $I:_ra_1\in \Mass_r(I)$, which
contradicts the maximality of $P$.
By induction there exists a sequence of elements of $a_1, a_2,\dots\notin I$
such that $I:_ra_1\subset I:_ra_2\subset \dots $. But this is a contradiction, because each $I:_ra_i$
is a right divisorial ideal and $S$ is a right Mori order. Thus $P$ is a maximal element of $\Gamma$.
\end{proof}

A right divisorial right ideal $I$ of a ring $R$ is called \textit{$\nu$-irreducible} 
if $I$ is not the intersection of two proper larger right divisorial right ideals.
Let $S$ be a right Mori order and $I$ a right divisorial right ideal. By Lemma~\ref{DCC},
the set of right divisorial right ideals properly containing  $I$ has minimal elements.
This set has a unique minimal element if and only if $I$ is $\nu$-irreducible.  
Following \cite{HL88} this unique minimal element is called the \textit{cover} of $I$.

The following is a generalization of \cite[Proposition 2.4]{HL88}.
\begin{proposition}\label{P2.8}
Let $S$ be a right Mori order is a simple Artinian ring $Q$ and $I$ a right divisorial right ideal of $S$. 
Then $I$ is a finite intersection of $\nu$-irreducible ideal of $S$. 
Moreover, each $\nu$-irreducible ideal has the form $aS:_rb$ for some regular elements $a, b\in S$.
\end{proposition}
\begin{proof}
Let $\Omega $ be the set of all minimal right divisorial right ideals properly containing $I$.
By Lemma~\ref{DCC}, the set $\Omega$ is not empty. Now if 
$|\Omega|=1$, then $I$ is {$\nu$-irreducible} and there is nothing to prove.
Hence we assume there exist $A, B \in\Omega$ with $A\neq B$. By minimality
of $A, B$ and the fact that the intersection of right divisorial right ideals is again right divisorial, 
we have $I=A\cap B$. Now if one of $A, B$, say $A$, is not $\nu$-irreducible, then 
there exist two properly larger right divisorial right ideals $A_1, A_2$ with 
$I\subset A\subset A_1, A_2$ and $I=A_1\cap A_2\cap B$. Continuing this procedure and using the fact that
$S$ has the ACC property on the set of right divisorial right ideals, we can find finitely many $\nu$-irreducible
ideals $J_1,\dots,J_n$ such that $I=J_i\cap \dots\cap J_n$.

To prove the second claim, let $I$ be $\nu$-irreducible and $J$ its cover. Since $I$ is right divisorial,
we have $I=\cap uS$, where $u\in U(Q)$ and $I\subseteq uS$. Since $I\subset J$, for at least one $u$ we have 
$J\nsubseteq uS\cap S$. From the minimality of $J$ and $\nu$-irreducibility of $I$ we conclude that $I=uS\cap S$.
Since $u\in U(Q)$, there exist regular elements $a, b\in S$ such that $u=b^{-1}a$.
Thus $I=b^{-1}aS\cap S=\{s\in S: s=b^{-1}ar ~\text{for some }~r\in S\}=
\{s\in S: bs=ar ~\text{for some }~r\in S\}=\{s\in S: bs\in aS\}=aS:_rb$. 
We can use the same procedure when none of $A$ and $B$ are $\nu$-irreducible.
\end{proof}

\section{examples of Mori orders}

Using an example from Cohen and Schofield \cite {CS85}, we will show that a right Mori order is not necessarily left Mori. 
 We refer the reader to \cite {CS85} and \cite {C85} for any undefined terminology. 
\begin{example}
Let $X=\{x_i : i\in \Z\}$, $Y=\{y\}$, $Z=\{z_i: i\in \Z\}$ and $Z_n=\{z_i: i\leq n\}$.
Put $R=\F\langle X, Y, Z: yz_i=z_{i-1}\rangle$ and $R_n=\F\langle X, Y, Z_n: yz_i=z_{i-1}\rangle$, 
the free $\F$-algebra on $X\cup Y\cup Z$ and $X\cup Y\cup Z_n$ respectively. Each $R_n$ is a \textit{free ideal ring (fir)},
i.e., all one sided ideals are free and of unique rank as left respectively right $R_n$-modules. By \cite [Theorem 10.3] {C85}, 
the ring $\F\langle y, z_i : yz_i=z_{i-1} \rangle$ is a right fir. We have   
$R=\F\langle y, z_i : yz_i=z_{i-1} \rangle \ast \F\langle X\rangle$, where $\ast$ denotes the coproduct.
$R$ is a right fir since it is the coproduct of two right firs. We recall that an $n\times n$ matrix over
a ring is called \textit{full} if it can not be written as the product of an $n\times (n-1)$ matrix and an $(n-1)\times n$ matrix.
Two matrices $A, B\in M_{n\times n}(S)$, where $S$ is a fir, are called \textit{totally coprime} if $A$ and $B$ have no common 
factors apart from units. Two sets of full matrices $\Delta$ and $\Lambda$ over a fir are called \textit{totally coprime}
if every element of $\Delta$ is totally coprime with every element of $\Lambda$.
Given $\Lambda=\{aI_{n\times n}: a\in Z\cup Y\}$ let $\Delta$ be the set of all full matrices over $R$ which are 
totally coprime to $\Lambda $. Let $\Delta_m$ be the subset of $\Delta$ with entries in $R_m$. Then we can 
construct the localization $R_{m\Delta_m}$ of $R_m$ by $\Delta_m$ which is defined as the ring obtained 
from $R_m$ by formally inverting all the matrices in $\Delta_m$. By \cite[Theorem 6.6] {C82}, $R_{m\Delta_m}$
is a simple principal ideal domain.
Let $U$ be the \textit{universal field of fractions} of $R$, that is $U$ is a field (not necessary commutative) 
with a homomorphism from $R$ to $U$ such that every full matrix over $R$ has an invertible image over $U$ \cite[cf. 7.2] {C85}.
Then all the ${R_m}_{\Delta_m}$ are subrings of $U$ as is $R_{\Delta}$. 
From $R_{1\Delta_1}\subset R_{2\Delta_2}\subset\dots$, $\cup R_{m\Delta_m}=R_{\Delta}$ and the fact that 
every $R_{m\,\Delta_m}$ is simple principal ideal domain, we conclude that $R_{\Delta}$ is a simple Bezout domain.
The ring $R$ is a right fir and the localization of a right fir is again right fir.
Hence $R_{\Delta}$ is a right principal ideal domain, so that $R_{\Delta}$ is right Mori.
Since $y, z_i$ are not units in $R_{\Delta}$ the ascending chain of left principal ideals 
$R_{\Delta}z_0\subset R_{\Delta}z_1\subset\dots $ cannot stabilize.
Therefore, $R_{\Delta}$ is not left Mori.
\end{example}

\medskip

A right $R$-module $M$ is called a \textit{generator} of the category of right $R$-modules
if $\sum_{f\in \Hom_R(M,R)} f(M)=R$. A \textit{progenerator} is a finitely generated projective right $R$-module $M$
which is a generator.

We recall the following well known theorem, see for example \cite[Theorem 1.9]{MMU97}. 
\begin{theorem}
Let $R$ and $S$ be rings. Then the following statements are equivalent:
\begin{itemize}
\item[(1)]
The categories of right $R$-modules and right $S$-modules are equivalent.
\item[(2)]
There exists a progenerator $P$ of the right $S$-module such that $R\cong \End_S(P)$. 
\item[(3)]
There exists an integer $n$ and idempotent $e\in M_{n\times n}(S)$ such that $R\cong eM_{n\times n}(S)e$
and $M_{n\times n}(S)eM_{n\times n}(S)=M_{n\times n}(S)$.
\end{itemize}
\end{theorem}

Two rings $R$ and $S$ are called \textit{Morita equivalent} if they satisfy the above conditions.
Properties of a module or a ring which are preserved under Morita equivalence are called \textit{Morita invariant}.

\begin{proposition}
Let $R$ be a right Mori order and let $S$ be a Morita equivalent to $R$. Then $S$ is right Mori. 
\end{proposition}
\begin{proof}
Let $F: \mathfrak{M}_R \rightarrow \mathfrak{M}_S$ and $G:  \mathfrak{M}_S\rightarrow \mathfrak{M}_R$ 
be naturally inverse category equivalences.
Let $Q=F(R_R)$ and $P=G(S_S)$. Then by \cite[(18.44) Proposition]{L99} the lattice of right ideals of $S$ (resp $R$)
is isomorphic to the lattice of submodules of $P_R$ (resp $Q_S$). This isomorphism maps a right ideal 
$I$ of $S$ (resp of $R$) to $I\otimes_SP\cong IP\subseteq P$ (resp  $I\otimes_RQ\cong IQ\subseteq Q$). Since 
$R$ is right Mori, there exists a finitely generated right ideal $J$ of $R$ such that $J\subseteq IP$
and $J^{\nu}=(IP)^{\nu}$. Thus $J\otimes _RQ\subseteq I\otimes_SP \otimes_R Q\cong I$ and $(J\otimes _RQ)^{\nu}=I^{\nu}$.
Since being finitely generated is a Morita invariant, the right ideal $J\otimes _RQ$ is finitely generated. Now by 
part (2) of Theorem \ref{chain}, $S$ is right Mori.
\end{proof}

\bibliographystyle{amsplain}

\end{document}